\documentclass[11pt,reqno]{amsart}

\usepackage{amssymb,latexsym,amsmath,amsfonts}
\usepackage{amsmath,amscd,amsthm,amssymb}

\numberwithin{equation}{section}

\def\whitebox{{\hbox{\hskip 1pt
 \vrule height 6pt depth 1.5pt
 \lower 1.5pt\vbox to 7.5pt{\hrule width
    3.2pt\vfill\hrule width 3.2pt}%
 \vrule height 6pt depth 1.5pt
 \hskip 1pt } }}
\def\qed{\ifhmode\allowbreak\else\nobreak\fi\hfill\quad\nobreak
     \whitebox\medbreak}
\newcommand{\pf}{\noindent{\bf Proof:}\ }

\newcommand{\ignore}[1]{}

\newcommand{\eproof}{\hfill$\Box$\vspace{4mm}}

\newcommand{\ga}{\gamma}
\newcommand{\F}{\mathbb{F}}
\newcommand{\Z}{\mathbb{Z}}

\newcommand{\Q}{\mathbb{Q}}
\newcommand{\la}{\langle}
\newcommand{\ra}{\rangle}

\newcommand{\Cc}{{\mathbb C}}

\newcommand{\Trace}{{\rm Tr}}

\newtheorem{thm}{Theorem}[section]

\newtheorem{example}[thm]{Example}
\numberwithin{equation}{section}

\begin{document}

\title[Pseudocyclic Association Schemes]
{Pseudocyclic and non-amorphic fusion schemes of the cyclotomic association schemes}

\author[Feng, Wu and Xiang]{Tao Feng$^*$, Fan Wu,  Qing Xiang$^\dagger$}

\dedicatory{Dedicated to Richard M. Wilson on the occasion of his 65th birthday}

\thanks{$^*$Supported in part by the Fundamental Research Funds for the central universities.}
      
\thanks{$^\dagger$Supported in part by NSF Grant DMS 1001557,  by the Overseas Cooperation Fund (grant 10928101) of China, and by Y. C. Tang disciplinary development fund, Zhejiang University.}

\address{Department of Mathematical Sciences, University of Delaware, Newark, DE 19716, USA.
{\bf Current Address:} Department of Mathematics, Zhejiang University, Hangzhou 310027, Zhejiang, China}
\email{pku.tfeng@yahoo.com.cn}

\address{Department of Mathematical Sciences, University of Delaware, Newark, DE 19716, USA}
\email{wufan@math.udel.edu}

\address{Department of Mathematical Sciences, University of Delaware, Newark, DE 19716, USA} \email{xiang@math.udel.edu}


\begin{abstract} 
We construct twelve infinite families of pseudocyclic and non-amorphic association schemes, in which each nontrivial relation is a strongly regular graph. Three of the twelve families generalize the counterexamples to A. V. Ivanov's conjecture by Ikuta and Munemasa \cite{ikutam}.
\end{abstract}

\maketitle

\section{Introduction}

This note is a sequel to \cite{FX}. We assume that the reader is familiar with the basic theory of association schemes as can be found in \cite{BI, BCN}. For background in strongly regular graphs, we refer the reader to \cite{bh, godsilr}. All association schemes considered in this paper are commutative and symmetric. Let $(X, \{R_i\}_{0\leq i\leq d})$ be an association scheme with $d$ classes. For $i\in \{0,1,\ldots ,d\}$, let $A_i$ be the adjacency matrix of the relation $R_i$, and let $E_0=\frac {1} {|X|}J, E_1, \ldots , E_d$ be the primitive idempotents of the Bose-Mesner algebra of the scheme $(X, \{R_i\}_{0\leq i\leq d})$, where $J$ is the all-one matrix of size $|X|\times |X|$. The basis transition matrix from $\{E_0,E_1,\ldots ,E_d\}$ to $\{A_0,A_1,\ldots ,A_d\}$ is denoted by $P=\left(p_j(i)\right)_{0\le i,j\le d}$, and usually called the {\it  first eigenmatrix} (or {\it character table}) of the scheme. Explicitly $P$ is  the $(d+1)\times (d+1)$ matrix with rows and columns indexed by $0,1,2,\ldots
,d$ such that
$$(A_0,A_1, \ldots ,A_d)=(E_0,E_1, \ldots ,E_d)P.$$ 
Let $k_i=p_i(0)$ and $m_i={\rm rank}(E_i)$. The $k_i$'s and $m_i$'s are called {\it valencies} and {\it multiplicities} of the scheme, respectively. We say that the scheme $(X, \{R_i\}_{0\leq i\leq d})$ is {\it pseudocyclic} if there exists an integer $t$ such that $m_i=t$ for all $i\in \{1,\ldots , d\}$. A classical example of pseudocyclic association schemes is the cyclotomic association scheme over a finite field, which we define below.

Let $q=p^f$, where $p$ is a prime and $f$ a positive integer. Let $\gamma$ be a fixed primitive element of $\F_q$ and $N|(q-1)$ with $N>1$. Let $C_0=\langle \gamma^N\rangle$, and $C_i=\gamma^i C_0$ for $1\leq i\leq N-1$. Assume that $-1\in C_0$. Define $R_0=\{(x,x) \mid  x\in \F_q\}$, and for $i\in \{1,2,\ldots ,N\}$, define $R_i=\{(x,y)\mid x,y\in \F_q, x-y\in C_{i-1}\}$. Then $(\F_q, \{R_i\}_{0\leq i\leq N})$ is an association scheme. We will call this scheme {\it the cyclotomic association scheme of class $N$ over $\F_q$}. The first eigenmatrix $P$ of the cyclotomic scheme of class $N$ is the following $(N+1)$ by $(N+1)$ matrix (with the rows of $P$ arranged in a certain way)
\begin{equation}\label{eigenmatrix}
P=\left(\begin{array}{cccccc}
1& \frac{N-1}{q}&\frac{N-1}{q}&\frac{N-1}{q}&\cdots &\frac{N-1}{q}\\
1&\eta_{_{N-1}}  &\eta_0   &\eta_1 & \cdots     &\eta_{_{N-2}}   \\
1&\eta_{_{N-2}}  &\eta_{_{N-1}} & \eta_0& \cdots  &\eta_{_{N-3}} \\
\vdots & & & & \\
1&\eta_0 &\eta_1 &\eta_2& \cdots  &\eta_{_{N-1}}\\
\end{array}\right)
\end{equation}
where the $\eta_i$'s are the cyclotomic periods  (or Gauss periods) of order $N$ defined by
$$\eta_i=\sum_{x\in C_i}\psi(x).$$ 
In the above defintion, $\psi$  is the additive character of $\F_q$ defined by
\begin{equation}\label{defaddchar}
\psi: \F_{q} \rightarrow \Cc^{*}, \quad \psi(x)=\xi_{p}^{\Trace(x)},
\end{equation}
where $\xi_p=e^{2\pi i/p}$ and ${\rm Tr}$ is the absolute trace from $\F_q$ to $\F_p$.

The following theorem gives combinatorial characterizations of pseudocyclic association schemes.

\begin{thm}\label{pseudocyc}
Let $(X, \{R_i\}_{0\leq i\leq d})$ be an association scheme, and for $x\in X$ and $1\leq i\leq d$, let $R_i(x)=\{y\mid (x,y)\in R_i\}$. Then the following are equivalent.
\begin{enumerate}
\item $(X, \{R_i\}_{0\leq i\leq d})$ is pseudocyclic.
\item For some constant $k$, we have $k_j=k$ and $\sum_{i=1}^{d}p_{ii}^j=k-1$, for $1\leq j\leq d$.
\item $(X, {\mathcal B})$ is a $2-(v,k,k-1)$ design, where ${\mathcal B}=\{R_i(x)\mid x\in X, 1\leq i\leq d\}$.
\end{enumerate}
\end{thm}
For a proof of this theorem, we refer the reader to \cite[p.~48]{BCN} and \cite[p.~84]{henkthesis}. Part (2) of the above theorem will be useful in Section 3.

Let $(X, \{R_i\}_{0\leq i\leq d})$ be an association scheme. 
For a partition $\Lambda_0:=\{0\}, \Lambda_1,\ldots ,\Lambda_{d'}$ of $\{0,1,\ldots ,d\}$, let $R_{\Lambda_i}=\cup_{k\in \Lambda_i}R_k$, 
for $0\leq i\leq d'$. If $(X, \{R_{\Lambda_i}\}_{0\leq i\leq d'})$ forms an association scheme, 
then we say that $(X, \{R_{\Lambda_i}\}_{0\leq i\leq d'})$ is a {\it fusion scheme} of the original scheme. 
If $(X, \{R_{\Lambda_i}\}_{0\leq i\leq d'})$ is an association scheme for every partition $\{\Lambda_i\}_{0\leq i\leq d'}$ of $\{0,1,2,\ldots ,d\}$ with $\Lambda_0=\{0\}$, 
then we call the original scheme $(X, \{R_i\}_{0\leq i\leq d})$ {\it amorphic}. For a recent survey on amorphic association schemes, we refer the reader to \cite{vanDamM}. Given a partition $\{\Lambda_i\}_{0\leq i\leq d'}$ of $\{0,1,2,\ldots ,d\}$ with $\Lambda_0=\{0\}$,  there is a simple criterion in terms of the first eigenmatrix $P$ of $(X, \{R_i\}_{0\leq i\leq d})$ for deciding whether $(X, \{R_{\Lambda_i}\}_{0\leq i\leq d'})$ forms an association scheme or not. We state this criterion below.\\

\noindent{\bf The Bannai-Muzychuk Criterion.} Let $P$ be the first eigenmatrix of an association scheme $(X, \{R_i\}_{0\leq i\leq d})$. Let $\Lambda_0:=\{0\}, \Lambda_1,\ldots ,\Lambda_{d'}$ be a partition of $\{0,1,\ldots ,d\}$. Then $(X, \{R_{\Lambda_i}\}_{0\leq i\leq d'})$ forms an association scheme if and only if there exists a partition $\{\Delta_i\}_{0\leq i\leq d'}$ of $\{0,1,2,\ldots ,d\}$ with $\Delta_0=\{0\}$ such that each $(\Delta_i, \Lambda_j)$-block of $P$ has a constant row sum. Moreover, the constant row sum of the $(\Delta_i, \Lambda_j)$-block is the $(i,j)$ entry of the first eigenmatrix of the fusion scheme. (For a proof of this criterion we refer the reader to \cite{BannaiSub, Muzthesis}.)

A. V. Ivanov conjectured in \cite{IP} that if each nontrivial relation in an association scheme is strongly regular, then the association scheme must be amorphic. 
This conjecture turned out to be false. A counterexample was given by Van Dam \cite{vanDamC1} in the case where the association scheme is imprimitive. 
Later on, Van Dam \cite{vanDamC2} also gave a counterexample  in the case where the association scheme is primitive.  
More counterexamples were given by Ikuta and Munemasa \cite{ikutam} in the primitive case. However it should be noted that there are only a few known counterexamples to Ivanov's conjecture in the primitive case (cf. \cite{ikutam}).

The purpose of this note is to generalize the counterexamples to Ivanov's conjecture by Ikuta and Munemasa \cite{ikutam} into infinite families. Along the way, we obtain many more infinite families of counterexamples to Ivanov's conjecture in the primitive case. The counterexamples we came up with are all pseudocyclic fusion schemes of the cyclotomic schemes. One of the main tools that we use is the theory of Gauss sums, which we review in the next section.

\section{Gauss sums}
Let $p$ be a prime, $f$ a positive integer, and $q=p^f$. Let $\xi_{p}=e^{2\pi i/p}$ and let $\psi$ be the additive character of $\F_{q}$ defined in (\ref{defaddchar}). Let
$$\chi:\F_{q}^* \rightarrow \Cc^{*} $$
be a character of $\F_{q}^{*}$. We define the {\it Gauss sum} by
$$ g(\chi)=\sum_{a \in \F_{q}^*} \chi(a)\psi(a).$$
Note that if $\chi_0$ is the trivial multiplicative character of
$\F_q$, then $g(\chi_0)=-1$. We are usually concerned with nontrivial Gauss sums $g(\chi)$, i.e., those with $\chi\neq \chi_0$.

While it is easy to show that the absoulte value of a nontrivial Gauss sum $g(\chi)$ is equal to $\sqrt{q}$, the explicit determination of Gauss sums 
is a difficult problem. However, there are a few cases where the Gauss sums $g(\chi)$ can be explicitly evaluated.  
The simplest case is the so-called {\it semi-primitive case}, where there exists an integer $j$ such that $p^j\equiv -1$ (mod $N$) ($N$ is the order of $\chi$ in ${\widehat \F_q^*}$, the character group of $\F_q^*$). Some authors \cite{BMW, bew} also refer to this case as uniform cyclotomy, or pure Gauss sums. We refer the reader to \cite[p.~364]{bew} for the precise evaluation of Gauss sums in this case. 

The next interesting case is the index 2 case, where $-1$ is not in the subgroup $\la p\ra$, the cyclic group generated by $p$,   
and $\la p\ra$ has index 2 in $(\Z/N\Z)^*$ (again here $N$ is the order of $\chi$ in ${\widehat \F_q^*}$). Many authors have investigated this case, including Baumert and Mykkeltveit \cite{bmy}, McEliece \cite{McE}, Langevin \cite{Lang}, Mbodj \cite{Mbo}, Meijer and Van de Vlugt \cite{mv}, and Yang and Xia  \cite{yx}. In the index 2 case, it can be shown that $N$ has at most two odd prime divisors. Assume that $N$ is odd, we have the following three possibilities in the index 2 case (see \cite{yx}): Below both $p_1$ and $p_2$ are primes.
\begin{enumerate}
\item $N=p_1^m$, $p_1\equiv 3$ (mod 4);
\item $N=p_1^mp_2^{n}$, $\{p_1\; ({\rm mod}\; 4), p_2\; ({\rm mod}\; 4)\}=\{1,3\}$, ${\rm ord}_{p_1^m}(p)=\phi(p_1^m)$, ${\rm ord}_{p_2^{n}}(p)=\phi(p_2^{n})$;
\item $N=p_1^mp_2^{n}$, $p_1\equiv 1,3$ (mod 4), ${\rm ord}_{p_1^m}(p)=\phi(p_1^m)$ and $p_2\equiv 3$ (mod 4), ${\rm ord}_{p_2^{n}}(p)=\phi(p_2^{n})/2$.
\end{enumerate}

We state below the results on evaluation of Gauss sums in Case (1) and (2)  from the above list.

\begin{thm}\label{lang}{\em (Langevin \cite{Lang})}
Let $N=p_1^m$, where $m$ is a positive integer, $p_1$ is a prime such that $p_1>3$ and $p_1\equiv 3$ (mod 4). Let $p$ be a prime such that $[(\Z/N\Z)^*:\la p\ra]=2$ (that is, $f:={\rm ord}_N(p)=\phi(N)/2$) and let $q=p^f$. Let $\chi$ be a  multiplicative character of order $N$ of $\F_q$, and $h$ be the class number of $\Q(\sqrt{-p_1})$. Then the Gauss sum $g(\chi)$ over $\F_q$ is determined up to complex conjugation by
$$g(\chi)=\frac{b+c\sqrt{-p_1}}{2}p^{h_0},$$
where 
\begin{enumerate}
\item $h_0=\frac {f-h}{2}$,
\item $b,c\not\equiv 0$ (mod $p$),
\item $b^2+p_1c^2=4p^{h}$,
\item $bp^{h_0}\equiv -2$ (mod $p_1$).
\end{enumerate}
\end{thm}

\begin{thm}\label{mbodj}{\em (Mbodj \cite{Mbo})}
Let $N=p_1^mp_2^{n}$, where $m,n$ are positive integers, $p_1$ and $p_2$ are prime such that $\{p_1\; ({\rm mod}\; 4), p_2\; ({\rm mod}\; 4)\}=\{1,3\}$, ${\rm ord}_{p_1^m}(p)=\phi(p_1^m)$, ${\rm ord}_{p_2^{n}}(p)=\phi(p_2^{n})$. Let $p$ be a prime such that $[(\Z/N\Z)^*:\la p\ra]=2$ (that is, $f:={\rm ord}_N(p)=\phi(N)/2$) and let $q=p^f$. Let $\chi$ be a multiplicative character of order $N$ of $\F_q$, and $h$ be the class number of $\Q(\sqrt{-p_1p_2})$. Then the Gauss sum $g(\chi)$ over $\F_q$ is determined up to complex conjugation by
$$g(\chi)=\frac{b+c\sqrt{-p_1p_2}}{2}p^{h_0},$$
where 
\begin{enumerate}
\item $h_0=\frac {f-h}{2}$,
\item $b,c\not\equiv 0$ (mod $p$),
\item $b^2+p_1p_2c^2=4p^{h}$,
\item $b\equiv 2p^{h/2}$ (mod $\ell$), here $\ell\in\{p_1,p_2\}$ is the prime congruent to 3 modulo 4.
\end{enumerate}
\end{thm}

\section{Pseudocyclic fusion schemes of the cyclotomic schemes}

Let $p$ be a prime, $f$ be a positive integer and $q=p^f$. Let $\gamma$ be a fixed primitive element of $\F_q$, and $N>1$ be an integer such that $N|(q-1)$. 
As we did in Section 1, let $C_0=\langle \gamma^N\rangle$ and $C_i=\gamma^iC_0$ for $1\leq i\leq N-1$. Assume that $-1\in C_0$. 
Define $R_0=\{(x,x) \mid  x\in \F_q\}$, and for $i\in \{1,2,\ldots ,N\}$, define $R_i=\{(x,y)\mid x,y\in \F_q, x-y\in C_{i-1}\}$. 
Then $(\F_q, \{R_i\}_{0\leq i\leq N})$ is the cyclotomic association scheme of class $N$ on $\F_q$. 
It was proven by Baumert, Mills and Ward \cite{BMW} that $(\F_q, \{R_i\}_{0\leq i\leq N})$ is amorphic if and only if $-1$ is congruent to a power of $p$ modulo $N$ 
(i.e., the so-called semi-primitive condition holds). See also \cite{BanMun} for a proof of this fact. Below we will show that even though in the index 2 case the cyclotomic association 
scheme $(\F_q, \{R_i\}_{0\leq i\leq N})$ is not amorphic, we can still have interesting fusion schemes of $(\F_q, \{R_i\}_{0\leq i\leq N})$.

\subsection{The index 2 case with $N=p_1^mp_2$}  

In this subsection, we assume that $N=p_1^mp_2$ ($m\geq 1$),  $p_1$, $p_2$ are primes such that $\{p_1\; ({\rm mod}\; 4), p_2\; ({\rm mod}\; 4)\}=\{1,3\}$, 
$p$ is a prime such that $\gcd(p, N)=1$, ${\rm ord}_{p_1^m}(p)=\phi(p_1^m)$ and ${\rm ord}_{p_2}(p)=\phi(p_2)$, and $f:={\rm ord}_N(p)=\phi(N)/2$. Let $q=p^f$, and as before let $C_0,C_1,\ldots ,C_{N-1}$ be the $N$-th cyclotomic classes of $\F_q$. Note that  here we have $-C_i=C_i$ for all $0\leq i\leq N-1$ since either $2N|(q-1)$ or $q$ is even. For convenience, we define $d:=p_1p_2$. For $0\leq k\leq d-1$, define
\begin{equation}
D_k=\bigcup_{i=0}^{p_1^{m-1}-1}C_{ip_2+kp_{1}^{m-1}}
\end{equation}
Note that $D_k=\gamma^{kp_1^{m-1}}D_0$ and $\{0\}, D_0, D_1, \ldots ,D_{d-1}$ form a partition of $\F_q$. Now define $R'_0=R_0$ and
\begin{equation}\label{defREL}
R'_k=\{(x,y)\mid x,y\in \F_q, x-y\in D_{k-1}\}.
\end{equation} 
We will show that $(\F_q, \{R'_k\}_{0\leq k\leq d})$ is a fusion scheme of $(\F_q,\{R_i\}_{0\leq i \leq N})$. The proof depends on the following evaluation of Gauss sums in the index 2 case, and results from \cite{FX}.

Let $\chi_1$ be the multiplicative character of order $p_1^m$ of $\F_q$ defined by $\chi_1(\ga)=\textup{exp}(\frac{2\pi i}{p_1^m})$, and let $\chi_2$ be the multiplicative character of order $p_2$ of $\F_q$ defined by $\chi_2(\ga)=\textup{exp}(\frac{2\pi i}{p_2})$. By Theorem~\ref{mbodj}, we have 
\begin{equation}\label{evalGauss}
g({\bar \chi_1}{\bar \chi_2})=\frac{b+c\sqrt{-p_1p_2}}{2}p^{h_0},
\end{equation}
where $h_0=\frac {f-h}{2}$ ($h$ is the class number of $\Q(\sqrt{-p_1p_2})$), $b,c\not\equiv 0$ (mod $p$), $b^2+p_1p_2c^2=4p^{h}$, and $b\equiv 2p^{h/2}$ (mod $\ell$), here $\ell\in \{p_1,p_2\}$ is the prime congruent to 3 modulo 4.

\begin{thm}\label{pseudoFusion1}
With the definition of $R'_k$ given in (\ref{defREL}), $(\F_q, \{R'_k\}_{0\leq k\leq d})$ is a pseudocyclic association scheme. 
\end{thm}

\pf We will first prove that $(\F_q, \{R'_k\}_{0\leq k\leq d})$ is an association scheme by using the Bannai-Muzychuk criterion discussed in Section 1. 

For each $a$, $0\le a\le N-1$, there exists a unique $i_{a}\in \{0, 1, ..., p_{1}^{m-1}-1\}$ such that $p_{1}^{m-1}\mid (a+p_{2}i_{a})$. It follows that there is a unique $j_{a}$, $0\leq j_a\leq p_1p_2-1$, such that $a\equiv -p_{2}i_{a}+p_{1}^{m-1}j_{a}$ (mod $N$). It is now easy to check that $-ip_2+jp_1^{m-1}$, $0\leq i\leq p_1^{m-1}-1$ and $0\leq j\leq p_1p_2-1$, form a complete set of residues modulo $N$.

The group of additive characters of $\F_q$ consists of $\psi_0$ and $\psi_{\gamma^a}$, $0\leq a\leq q-2$,  where $\psi_0$ is the trivial character and $\psi_{\gamma^a}$ is defined by
\begin{equation}
\psi_{\gamma^a}: \F_{q} \rightarrow \Cc^{*}, \quad \psi_{\gamma^a}(x)=\xi_{p}^{\Trace(\gamma^ax)}.
\end{equation}
We usually write $\psi_1$ simply as $\psi$. The character values of $D_0$ were computed in the proof of Theorem 5.1 \cite{FX}. Since $D_k$ is a (multiplicative) translate of $D_0$, we know the character values of $D_k$ as well. Explicitly, for each $a$, $0\leq a\leq N-1$, write 
$$a\equiv -p_{2}i_{a}+p_{1}^{m-1}j_{a}\; ({\rm mod}\; N),$$
with $0\leq i_a\leq p_1^{m-1}-1$ and $0\leq j_a\leq p_1p_2-1$. For convenience we introduce the Kronecker delta $\delta_{a,p_1}$, which equals 1 if $p_1|a$, 0 otherwise. Also we define $\delta_{a,p_2}$ by setting it equal to 1 if $p_2|a$, 0 otherwise. By the results in \cite{FX}, we have $$\psi_{\gamma^a}(D_k)=\psi(\gamma^{a+p_1^{m-1}k}D_0)=\frac{1}{N}T_{a+p_1^{m-1}k},$$ 
where
\begin{align*}
T_{a+p_1^{m-1}k}&=-p_1^{m-1}-(-1)^{\frac{p_1-1}{2}}p_1^{m-1}p_2\sqrt{q}\delta_{a+p_1^{m-1}k,p_2}-(-1)^{\frac{p_2-1}{2}}p_1^m\sqrt{q}\delta_{j_a+k,p_1}\\
&\quad+\frac{b}{2}p^{h_0}p_1^{m-1}(p_1\delta_{j_a+k,p_1}-1)(p_2\delta_{a+p_1^{m-1}k,p_2}-1)\\
&\quad -\bigg(\frac{a+p_1^{m-1}k}{p_2}\bigg)\bigg(\frac{j_a+k}{p_1}\bigg)\frac{c}{2}p^{h_0}p_1^mp_2
\end{align*}
In the above formula, $b,c$ are given by (\ref{evalGauss}), $(\frac{.}{p_2})$ and $(\frac{.}{p_1})$ are Legendre symbols. 
Observe that $a+p_1^{m-1}k\equiv -p_{2}i_{a}+p_{1}^{m-1}(j_{a}+k)$ (mod $N$). So $\delta_{a+p_1^{m-1}k, p_{2}}=\delta_{j_{a}+k,p_2}$, 
and $\bigg(\frac{a+p_1^{m-1}k}{p_{2}}\bigg)=\bigg(\frac{p_1}{p_2}\bigg)^{m-1}\bigg(\frac{j_a+k}{p_2}\bigg)$. Therefore, $\psi_{\gamma^a}(D_k)$ is independent of $i_a$.\\

In order to apply the Bannai-Muzychuk criterion, we define the following partition of $\{\psi_{\gamma^a}\mid a\in \Z/N\Z\}$. For each $j$, $0\leq j\leq d-1$, define
$$\Delta_{j+1}=\{\psi_{\gamma^{-p_2i + p_1^{m-1} j}}\mid  0\leq i\leq p_1^{m-1}-1\},$$ 
and $\Delta_{0}=\{\psi_0\}$. Clearly $\Delta_{0}, \Delta_{1}, \ldots ,\Delta_{d}$ form a partition of $\{\psi_{\gamma^a}\mid a\in \Z/N\Z\}$.   For each $0\leq k\leq d-1$, since $\psi_{\gamma^a}(D_k)$ is independent of $i_a$ (here $a\equiv -p_{2}i_{a}+p_{1}^{m-1}j_{a}\; ({\rm mod}\; N)$),  we see that $\psi_{\gamma^a}(D_k)$ is a constant for those $a$ in the same subset of the above partition. By the Bannai-Muzychuk criterion (with $\Lambda_0=\{0\}$, $\Lambda_{j+1}=\{1+ip_2+p_1^{m-1}j\mid 0\leq i\leq p_1^{m-1}-1\}$, $0\leq j\leq d-1$), we see that $(\F_q, \{R'_{0}$, $R'_1,\ldots , R'_{d}\})$ is an association scheme.\\

Next we show that the association scheme $(\F_q, \{R'_k\}_{0\leq k\leq d})$ is pseudocyclic. To this end, we show that the following group ring equation holds in $\Z[(\F_q,+)]$.\\

\noindent{\it Claim:} $\sum_{k=0}^{d-1}D_k^2=(q-1)\cdot 0_{\F_q}+\frac{q-1}{p_1p_2}(\F_q-0_{\F_q})$, where $0_{\F_q}$ is the zero element in $\F_q$. \\

For any $a$, $0\leq a\leq N-1$, we write $a\equiv -i_ap_2+j_ap_1^{m-1} \pmod{N}$ with $i_a\in\{0,1,\ldots,p_1^{m-1}-1\}$ and $j_a\in\{0,1,2,\ldots ,d-1\}$. Since $\psi_{\gamma^a}(D_k)$ is independent of $i_a$, we may assume that $i_a=0$. We now compute
\begin{align*}
 \sum_{k=0}^{d-1}\psi_{\gamma^a}(D_k)^2=\frac{1}{N^2}\sum_{k=0}^{d-1}T_{p_1^{m-1}(j_a+k)}^2=\frac{1}{N^2}\sum_{k=0}^{d-1}T_{kp_1^{m-1}}^2
\end{align*}
Since the last expression above is independent of $a$, we see that $\sum_{k=0}^{d-1}\psi_{\gamma^a}(D_k)^2$ are equal to the same constant for all $0\leq a\leq N-1$. Since each $D_k$ is a union of some $N$-th cyclotomic classes, it follows that $\sum_{k=0}^{d-1}\psi_{\gamma^a}(D_k)^2$ are equal to the same constant for all $0\leq a\leq q-2$. Therefore, by the inversion formula, we have 
\[\sum_{k=0}^{d-1}D_k^2=(n-\lambda)\cdot 0_{\F_q}+\lambda\F_q,\]
for some integers $n,\lambda$. Now applying the principal character to both sides, and computing the coefficients of $0_{\F_q}$ on both sides, we have 
\begin{align*}
 n&=p_1p_2\cdot\frac{q-1}{p_1p_2},\\
 n+(q-1)\lambda&=d\cdot\bigg(\frac{q-1}{p_1p_2}\bigg)^2.
\end{align*}
It follows that $n=q-1$, and $\lambda=\frac{q-1}{p_1p_2}-1$. The claim is now established. A direct consequence is that $\sum_{i=0}^{d-1}p_{i,i}^j=\frac{q-1}{N}-1$, for all $j$, where $p_{i,i}^j$ are the intersection parameters given by  $D_i^2=\sum_{j=0}^{d-1}p_{i,i}^jD_j+p_{i,i}^{0}\cdot 0_{\F_q}$. By Part (2) of Theorem~\ref{pseudocyc}, the association scheme $(\F_q, \{R'_k\}_{0\leq k\leq d})$ is pseudocyclic. The proof is complete.
\eproof

In order to obtain counterexamples to Ivanov's conjecture, we need to have each $R'_k$ ($1\leq k\leq d$) in Theorem~\ref{pseudoFusion1} to be strongly regular. 
Note that $R'_k$ is just the Cayley graph ${\rm Cay}(\F_q, D_{k-1})$, and ${\rm Cay}(\F_q, D_{k-1})\cong {\rm Cay}(\F_q, D_0)$ for all $1\leq k\leq d$ since $D_{k-1}=\gamma^{(k-1)p_1^{m-1}}D_0$. 
It follows that if ${\rm Cay}(\F_q, D_0)$ is strongly regular, then all $R'_k$, $1\leq k\leq d$, are strongly regular. In \cite{FX}, we obtained necessary and sufficient conditions for ${\rm Cay}(\F_q, D_0)$ to be strongly regular, which we quote below.

\begin{thm}{\rm (Corollary 5.2 in \cite{FX})}\label{srg2}
With $b,c,h$ given in (\ref{evalGauss}), ${\rm Cay}(\F_q,D_0)$ is a strongly regular graph if and only if $b,c\in \{1,-1\}$, $h$ is even and $p_1=2p^{h/2}+(-1)^{\frac{p_1-1}{2}}b$, $p_2=2p^{h/2}-(-1)^{\frac{p_1-1}{2}}b$. 
\end{thm}

In \cite{FX}, we used a computer to search for $p, p_1, p_2$ satisfying the conditions in Theorem~\ref{srg2}. We found six infinite families of strongly regular graphs in this way. By the discussion preceding Theorem~\ref{srg2}, and since the parameters of each of the six examples of srg are neither Latin square type nor negative Latin square type, each of the six families of srg gives rise to an infinite class of counterexamples to Ivanov's conjecture. Below we list the parameters of these examples. For the detailed reasons why we have strongly regular graphs, we refer the reader to \cite{FX}.

\begin{example}\label{21st}
Let $p=2$, $q=2^{4\cdot 3^{m-1}}$, $p_1=3$, $p_2=5$, $N=3^m\cdot 5$, with $m\geq 1$. Then we have a 15-class pseudocyclic fusion scheme $(\F_q, \{R'_k\}_{0\leq k\leq 15})$ in which each relation $R'_k$, $1\leq k\leq 15$, is strongly regular. 
\end{example}

We remark that when $m=2$, Example~\ref{21st} is the same as Example 1 in \cite{ikutam}.

\begin{example}\label{22st}
Let $p=2$, $q=2^{4\cdot 5^{m-1}}$, $p_1=5$, $p_2=3$, $N=5^m\cdot 3$, with $m\geq 1$. Then we have a 15-class pseudocyclic fusion scheme $(\F_q, \{R'_k\}_{0\leq k\leq 15})$ in which each relation $R'_k$, $1\leq k\leq 15$, is strongly regular. 
\end{example}

We remark that when $m=2$, Example~\ref{22st} is the same as Example 2 in \cite{ikutam}.

\begin{example}\label{23st}
Let $p=3$, $q=3^{12\cdot 5^{m-1}}$, $p_1=5$, $p_2=7$, $N=5^m\cdot 7$, with $m\geq 1$. Then we have a $35$-class pseudocyclic fusion scheme $(\F_q, \{R'_k\}_{0\leq k\leq 35})$ in which each relation $R'_k$, $1\leq k\leq 35$, is strongly regular. 
\end{example}

\begin{example}\label{24st}
Let $p=3$, $q=3^{12\cdot 5^{m-1}}$, $p_1=7$, $p_2=5$, $N=7^m\cdot 5$, with $m\geq 1$. Then we have a $35$-class pseudocyclic fusion scheme $(\F_q, \{R'_k\}_{0\leq k\leq 35})$ in which each relation $R'_k$, $1\leq k\leq 35$, is strongly regular. 
\end{example}

\begin{example}\label{25st}
Let $p=3$, $q=3^{144\cdot 17^{m-1}}$, $p_1=17$, $p_2=19$, $N=17^m\cdot 19$, with $m\geq 1$. Then we have a $323$-class pseudocyclic fusion scheme $(\F_q, \{R'_k\}_{0\leq k\leq 323})$ in which each relation $R'_k$, $1\leq k\leq 323$, is strongly regular. 
\end{example}

\begin{example}\label{26st}
Let $p=3$, $q=3^{144\cdot 19^{m-1}}$, $p_1=19$, $p_2=17$, $N=19^m\cdot 17$, with $m\geq 1$. Then we have a $323$-class pseudocyclic fusion scheme $(\F_q, \{R'_k\}_{0\leq k\leq 323})$ in which each relation $R'_k$, $1\leq k\leq 323$, is strongly regular.
\end{example}

We remark that by using Corollary 3.2 in \cite{ikutam}, one can further obtain 3-class fusion schemes of the above pseudocyclic association schemes, in which two relations are strongly regular graphs, while the third relation is not (see the character table of these 3-class fusion schemes in the statement of Corollary 3.2 of \cite{ikutam}).

\subsection{The index 2 case with $N=p_1^m$} 

In this subsection, we assume that $N=p_1^m$  (here $m\geq 1$,  $p_1>3$ is a prime such that $p_1\equiv 3$ (mod 4)), $p$ is a prime such that $\gcd(N,p)=1$, and  $f:={\rm ord}_{N}(p)=\phi(N)/2$. Let $q=p^f$, and as before let $C_0,C_1,\ldots ,C_{N-1}$ be the $N$-th cyclotomic classes of $\F_q$. Note that $-C_i=C_i$ for all $0\leq i\leq N-1$ since either $2N|(q-1)$ or $q$ is even. For $0\leq k\leq p_1-1$, define
\begin{equation}
D_k=\bigcup_{i=0}^{p_1^{m-1}-1}C_{i+kp_1^{m-1}}
\end{equation}
Note that $D_k=\gamma^{kp_1^{m-1}}D_0$ and $\{0\}, D_0, D_1, \ldots ,D_{p_1-1}$ form a partition of $\F_q$. Now define $R'_0=R_0$ and
\begin{equation}\label{defREL2}
R'_k=\{(x,y)\mid x,y\in \F_q, x-y\in D_{k-1}\}.
\end{equation} 
We will show that $(\F_q, \{R'_k\}_{0\leq k\leq p_1})$ is a fusion scheme of $(\F_q,\{R_i\}_{0\leq i \leq N})$. The proof depends on the following evaluation of Gauss sums in the index 2 case, and results from \cite{FX}.

Let $\chi$ be the multiplicative character of $\F_q$ defined by $\chi(\ga)=\text{exp}(\frac{2\pi i}{N})$. By Theorem~\ref{lang}, we have
\begin{equation}\label{evalGauss2}
g({\bar \chi})=\frac{b+c\sqrt{-p_1}}{2}p^{h_0},\;\; b,c\not\equiv 0\pmod{p},
\end{equation}
where $h_0=\frac {f-h}{2}$ and $h$ is the class number of $\Q(\sqrt{-p_1})$, $b^2+p_1c^2=4p^{h}$, and $bp^{h_0}\equiv -2\pmod{p_1}$. 

\begin{thm}\label{pseudoFusion2}
With the definition of $R'_k$ given in (\ref{defREL2}), $(\F_q, \{R'_k\}_{0\leq k\leq p_1})$ is a pseudocyclic association scheme. 
\end{thm}

\pf The proof is similar to that of Theorem~\ref{pseudoFusion1}. For each $a$, $0\leq a\leq N-1$, there is a unique $i_a\in\{0,1,\ldots,p_1^{m-1}-1\}$, such that $p_1^{m-1}|(a+i_a)$. It follows that there is a unique $j_a, 0\leq j_a\leq p_1-1$, such that $a\equiv -i_a+p_1^{m-1}j_a\; ({\rm mod}\; N).$ It is now easy to check that $-i+jp_1^{m-1}$, $0\leq i\leq p_1^{m-1}-1$ and $0\leq j\leq p_1-1$, form a complete set of residues modulo $N$.

The group of additive characters of $\F_q$ consists of $\psi_0$ and $\psi_{\gamma^a}$, $0\leq a\leq q-2$. The character values of $D_0$ were computed in the proof of Theorem 4.1 \cite{FX}.  Since $D_k$ is a (multiplicative) translate of $D_0$, we know the character values of $D_k$ as well. Explicitly, for each $a$, $0\leq a\leq N-1$, write 
$$a\equiv -i_a+p_1^{m-1}j_a\; ({\rm mod}\; N),$$
with $0\leq i_a\leq p_1^{m-1}-1$ and $0\leq j_a\leq p_1-1$. For convenience, we also introduce the Kronecker delta $\delta_{j_a}$, which equals $1$ if $p_1|j_a$, and $0$ otherwise. By the results in \cite{FX}, we have 
$$\psi_{\gamma^a}(D_k)=\psi(\gamma^{a+kp_1^{m-1}}D_0)=\frac{1}{N}T_{a+kp_1^{m-1}},$$ 
where 
$$T_{a+kp_1^{m-1}}=-p_1^{m-1}+\frac{p^{h_0}p_1^{m-1}b}{2}(p_1\delta_{j_a+k}-1)-\frac{p^{h_0}p_1^mc}{2}\bigg(\frac{j_a+k}{p_1}\bigg).$$ 
In the above formula, $b,c$ are given in (\ref{evalGauss2}), and $(\frac{.}{p_1})$ is the Legendre symbol. It is important to note that $\psi_{\gamma^a}(D_k)$ is independent of $i_a$. \\

We define the following partition of $\{\psi_{\gamma^a}\mid a\in \Z/N\Z\}$. For each $j$, $0\leq j\leq p_1-1$, we define
$$\Delta_{j+1}=\{\psi_{\gamma^{-i + p_1^{m-1} j}}\mid  0\leq i\leq p_1^{m-1}-1\},$$ 
and $\Delta_{0}=\{\psi_0\}$. Then clearly $\Delta_{0}, \Delta_{1}, \ldots ,\Delta_{p_1}$ form a partition of $\{\psi_{\gamma^a}\mid a\in \Z/N\Z\}$.   For each $0\leq k\leq p_1-1$, since $\psi_{\gamma^a}(D_k)$ is independent of $i_a$ (here $a\equiv -i_{a}+p_{1}^{m-1}j_{a}\; ({\rm mod}\; N)$),  we see that $\psi_{\gamma^a}(D_k)$ is a constant for those $a$ in the same subset of the above partition. By the Bannai-Muzychuk criterion (with $\Lambda_0=\{0\}$, $\Lambda_{j+1}=\{1+i+p_1^{m-1}j\mid 0\leq i\leq p_1^{m-1}-1\}$, $0\leq j\leq p_1-1$), we see that $(\F_q, \{R'_{0}$, $R'_1\ldots , R'_{p_1}\})$ is an association scheme.

Similarly we can show that the following group ring equation holds in $\Z[(\F_q,+)]$:$$\sum_{k=0}^{p_1-1}D_k^2=(q-1)\cdot 0_{\F_q}+\frac{q-1}{p_1}(\F_q-0_{\F_q}),$$ 
from which the pseudocyclicity of the scheme $(\F_q, \{R'_{0}$, $R'_1,\ldots , R'_{p_1}\})$ follows. We omit the details of the proof of the above group ring equation. The proof is now complete.
\eproof

In order to obtain counterexamples to Ivanov's conjecture, we need to have each $R'_k$ ($1\leq k\leq p_1$) in Theorem~\ref{pseudoFusion2} to be strongly regular. 
Note that $R'_k$ is just the Cayley graph ${\rm Cay}(\F_q, D_{k-1})$, and ${\rm Cay}(\F_q, D_{k-1})\cong {\rm Cay}(\F_q, D_0)$ for all $1\leq k\leq p_1$ since $D_{k-1}=\gamma^{(k-1)p_1^{m-1}}D_0$. 
Again it follows that if ${\rm Cay}(\F_q, D_0)$ is strongly regular, then all $R'_k$, $1\leq k\leq p_1$, are strongly regular. 
In \cite{FX}, we obtained necessary and sufficient conditions for ${\rm Cay}(\F_q, D_0)$ to be strongly regular, which we quote below.

\begin{thm}{\rm (Corollary 4.2 in \cite{FX})}\label{srg1}
With $b,c$ given in (\ref{evalGauss2}), ${\rm Cay}(\F_q,D)$ is a strongly regular graph if and only if $b,c\in \{1,-1\}$.
\end{thm}

In \cite{FX}, we used a computer to search for $p, p_1$ satisfying the conditions in Theorem~\ref{srg1}. We found six infinite families of strongly regular graphs in this way. 
By the discussion preceding Theorem~\ref{srg1}, each of the six examples of srg gives rise to a class of infinitely many counterexamples to Ivanov's conjecture. Below we list the parameters of these examples. For the detailed reasons why we have strongly regular graphs, we refere the reader to \cite{FX}.

\begin{example}\label{1st}
Let $p=2$, $q=2^{3\cdot 7^{m-1}}$, $p_1=7$, $N=p_1^m$, $m\geq 1$ is an integer. Then we have a $7$-class pseudocyclic fusion scheme $(\F_q, \{R'_k\}_{0\leq k\leq 7})$ in which each relation $R'_k$, $1\leq k\leq 7$, is strongly regular. 
\end{example}

We remark that when $m=2$, Example~\ref{1st} is the same as Example 3 of \cite{ikutam}.

\begin{example}\label{2nd}
Let $p=3$, $q=3^{53\cdot 107^{m-1}}$, $p_1=107$, $N=p_1^m$, $m\geq 1$ is an integer. Then we have a $107$-class pseudocyclic fusion scheme $(\F_q, \{R'_k\}_{0\leq k\leq 107})$ in which each relation $R'_k$, $1\leq k\leq 107$, is strongly regular. 
\end{example}

\begin{example}\label{3rd}
Let $p=5$, $q=5^{9\cdot 19^{m-1}}$, $p_1=19$, $N=p_1^m$,  $m\geq 1$ is an integer. Then we have a $19$-class pseudocyclic fusion scheme $(\F_q, \{R'_k\}_{0\leq k\leq 19})$ in which each relation $R'_k$, $1\leq k\leq 19$, is strongly regular. 
\end{example}

\begin{example}\label{4th}
Let $p=5$, $q=5^{249\cdot 499^{m-1}}$, $p_1=499$, $N=p_1^m$, $m\geq 1$ is an integer. Then we have a $499$-class pseudocyclic fusion scheme $(\F_q, \{R'_k\}_{0\leq k\leq 499})$ in which each relation $R'_k$, $1\leq k\leq 499$, is strongly regular. 
\end{example}

\begin{example}\label{5th}
Let $p=17$, $q=17^{33\cdot 67^{m-1}}$, $p_1=67$,  $N=p_1^m$, $m\geq 1$ is an integer. Then we have a $67$-class pseudocyclic fusion scheme $(\F_q, \{R'_k\}_{0\leq k\leq 67})$ in which each relation $R'_k$, $1\leq k\leq 67$, is strongly regular. 
\end{example}

\begin{example}\label{6th}
Let $p=41$, $q=41^{81\cdot 163^{m-1}}$, $p_1=163$, $N=p_1^m$, $m\geq 1$ is an integer. Then we have a $163$-class pseudocyclic fusion scheme $(\F_q, \{R'_k\}_{0\leq k\leq 163})$ in which each relation $R'_k$, $1\leq k\leq 163$, is strongly regular. 
\end{example}

Again we remark that by using Corollary 3.2 in \cite{ikutam}, one can further obtain 3-class fusion schemes of the above pseudocyclic association schemes, in which two relations are strongly regular graphs, while the third relation is not.

\end{document}